\documentclass{amsart}
\usepackage{amssymb, amsmath, amsthm}

\title{Thomae formula for general cyclic covers of $\mathbb{CP}^1$ }
\author[Kopeliovich]{by Yaacov Kopeliovich}
\date{23 September 2010 }
\address{MEAG NY, 540 Madison Avenue sixth floor New York NY 10022 and Department of Mathematics Bronx Community College 2155 University Avenue Bronx NY 10453}
\email{ykopeliovich@yahoo.com}
\newtheorem{thm}{Theorem}[section]
\newtheorem{cor}[thm]{Corollary}
\newtheorem{lem}[thm]{Lemma}

\newtheorem{prop}[thm]{Proposition}
\theoremstyle{definition}
\newtheorem{defn}[thm]{Definition}

\theoremstyle{remark}
\newtheorem{rmk}[thm]{Remark}

\begin{document}
\begin{abstract}
Let $X$ be a general cyclic cover of  $\mathbb{CP}^{1}$ ramified at $m$ points, $\lambda_1...\lambda_m.$ we define a class of non positive divisors on $X$ of degree $g-1$ supported in the pre images of the branch points on $X$, such that the  Riemann theta function doesn't vanish on their image in $J(X).$  We generalize the results of [BR],[Na] and [EG] and prove that up to a certain determinant of the non standard periods of $X$, the value of the Riemann theta function at these divisors raised to a high enough power is a polynomial in the branch point of the curve $X.$ Our approach is based on a refinement of Accola's results for 3 cyclic sheeted cover [Ac1] and a generalization of Nakayashiki's approach explained in [Na] for general cyclic covers.
\end{abstract}
\maketitle
\section{Introduction}
Let $X$ be an algebraic curve given by the equation :
$$y^N=\prod_{i=}^m \left(x-\lambda_i\right)^{R_i}$$ such that $\sum_{i=1}^m R_i=0\bmod N$ and $(R_i,N)=1.$ Let $\phi:X\mapsto \mathbb{CP}^{1}$ be a map defined by: $\phi(x,y)=x.$  This map is $N$ to $1$ from $X$ to $\mathbb{CP}^{1}.$ $\lambda_i$ are the ramification points. Let $P_i=\phi^{-1}(\lambda_i).$
Choose a base point $z_0,$ a normalized homology basis $a_1...a_g,b_1...b_g$ and a normalized holomorphic differentials $v_1,...v_g$ to define the Jacobian of $J(X)$ and a standard map $u:X\mapsto J(X).$ Let $K_{z_0}$ be the Riemann constant and $\tau$ is the period matrix associated with the homology basis and the differentials selected above. We prove the following theorem:
\begin{thm}
Let $r$ be a total ramification of $\phi:X\mapsto \mathbb{CP}^{1}.$ Select an integer vector $\beta=\left(\beta_1...\beta_m\right)$ such that:
\begin{enumerate}
  \item $0\leq \beta_i\leq N-1$
  \item for $0\leq k\leq N-1,$ $\sum_{i=1}^m \overline{\beta_i+kR_i}=\frac{r}{2}$
\end{enumerate}
and  $\overline{j}$ denotes the smallest positive integer $j_0$ such that $j\bmod N=j_0 \bmod N.$

Let $K_{z_0}$ be the Riemann's constant and let $\theta(z,\tau)$ be the Riemann theta function.
Then $$\theta\left[u(\sum_{i=1}^m \beta_iP_i)+ K_{z_0}-u(\sum_{j=1}^N\infty_j)\right](0,\tau)\neq 0$$ and there exists a complex number $\alpha$ not depending on $\tau$ such that:
$$
\theta\left[u(\sum_{j=1}^m \beta_iP_i)+K_{z_0}-u(\sum_{i=1}^N\infty_i)\right]=
\alpha\sqrt{\det
C}\times{\prod_{i,j=1..m,i\neq j}(\lambda_i-\lambda_j)}^{\beta_{ij}+\frac{\gamma_{ij}}{2}}
$$
where  $\det C$ is a certain determinant of the $g\times g$ matrix of non normalized holomorphic differentials evaluated at $a_i, 1\leq i\leq g$ and
        $$\gamma_{ij}=\sum_{w=0}^{N-1}\left\{wR_i/N\right\}\left\{wR_j/N\right\},$$
         $$\beta_{ij}=\sum_{k=0}^{N-1}\left(\left\{(\beta_i+kR_i)/N\right\}-\frac{N-1}{2N}\right)\times\left(\left\{(\beta_j+kR_j)/N\right\}-\frac{N-1}{2N}\right)$$
and $\{\alpha\}$ is the fractional part of $\alpha.$
\end{thm}
the theorem is a generalization of the work started by [BR],[Na] and [EG]. Using methods from String and Quantum field theory Bershadsky and Radul generalized Thomae formula for hyper-elliptic covers to a non singular covers of the sphere i.e. when $R_i=1$ and the number of branch points $m$ is a multiple of $N.$ [Na] gave a more rigorous proof for the formula suggested by [BR] while [EG] modified Nakayashiki's method and treated a special singular case. In this note we follow the  approach of [Na] to prove the formula we stated above. First we show that the Riemann theta function doesn't vanish on the images of the divisors we defined above. Then we modify [Na] to the case when $(R_i,N)=1$ and $\sum_{i=1}^m R_i=0 \bmod N.$ As far as we know it provides Thomae formula for widest class of cyclic covers of the Sphere.

The main idea of [Na] is to produce an integrable differential equation that describes the variation of the logarithm of the theta function with respect to the branch points. [Na] constructs certain analytic quantities of the Riemann surfaces locally ( as algebraic expressions supported by local coordinates around the branch points) and compares them to the global expression as derived in [Fa]. Equating the expansions of the global and the local constructions produces the result. We carry this program below. It turns out that the general case of cyclic covers we handle doesn't differ much from the case considered by [Na] and [BR].

It is interesting to look for Thomae formulas since they should be useful in mathematical physics, representation and number theory. [Na1] applied the formulas he found in a non singular case to investigate solutions of KZ equations. [EG] relied on the formula the Thomae formula they found to solve the Riemann Hilbert problem for special class of cyclic covers. On the other hand these formulas might be of interest in representation theory, since Symmetric groups (or products of Symmetric groups), act on polynomials on the RHS of the formula. Consequently powers of the  theta functions realize representations of Symmetric groups ( or  their products). One can then use the machinery of representation theory to derive a basis for the theta functions and enhance our  understanding of their modular properties and perhaps characterize the periods coming from cyclic of covers of $\mathbb{CP}^1$ . For this approach, when the cyclic cover is of degree, see [Ko]. In number theory Thomae formulas can probably be used to improve existing algorithms for counting points on cyclic covers of the projective line above finite fields. Mestre applied Thomae formula and duplication formulas of the theta functions count points of these curves above finite fields. It is plausible that his approach can  be generalized to the cyclic cover case as well using the generalized Thomae formula.
Finally one can't ignore the inherent mystery of the formulas where the LHS is a highly transcendental object while a RHS is a product of differences of points.

The approach presented here isn't the only one to look for  these formulas. In a series of papers Hershel Farkas and his collaborators ([EF],[EiF]) reproved Thomae's original result and used classical approach of Riemann to write the branch points as of cyclic covers as ratios of theta functions. Consequently they were able to get the $\beta$ part of the formula in certain cases. In a book currently written with his student Zemel, Farkas [FZ] generalizes his work avoiding the variational approach used in this note. Lastly [KT] derived a similar result for the case $p=3.$ More precisely [KT] assumed Nakayashiki's result for a non singular cover cyclic cover of degree $3$ and analyzed degenerations of the homology basis when ramification points coalesce. Their approach enables them to calculate the constant $\alpha$ explicitly. In subsequent notes we will attempt to understand  their method and reprove results obtained in this note.

\thanks{This work was partially done during a visit to the TAMU math
department and the author thanks the department for the invitation
and kind hospitality. He thanks Hershel Farkas for very useful discussions on Thomae formula and he especially thanks
 Mike Fried for pointing out an inaccuracy in the draft version of this paper and generously sharing his ideas for a possible future research. Many thanks to the referee for carefully checking the manuscript
and pointing out numerous corrections in the mathematical and
stylistic content of this note.
}

\section{Non positive divisors on Riemann surface}
Let $X$ be a Riemann surface and assume that $D=\sum d_iz_i$ is a divisor ( not necessarily positive) on it.
\begin{defn}
$H^0\left(X,\mathcal{O}\left(D\right)\right)$ is the collection of functions $f:X\mapsto \mathbb{CP}^1$ on $X$ such that $div(f)\geq D.$ In the [FK] notation this is the space $\mathfrak{R}\left(D\right).$ \\Let $r(D)= dim H^0\left(X,\mathcal{O}\left(D\right)\right).$
\end{defn}
We seek conditions when $\exists E$ a divisor on $X$ such that $E=\sum e_ix_i,$   $e_i\geq 0$ and $D \equiv E.$ Assume that $D$ is not a positive divisor (otherwise you can set $E=D$) Then if $E$ is positive and equivalent to $D$ there exists a non constant function $f$ such that $div(f)=E/D.$ Therefore $f \in H^0\left(X,\mathcal{O}\left(-D\right)\right).$ (That is $f$ is a function such that $div(f) \geq D.)$ Conclude that $r(-D)>0.$ We showed the following: \begin{lem}
Let $D$ be a non positive divisor. Then if there is $E$ a positive divisor such that $E\equiv D$ then $r(-D)>0.$
\end{lem}
Note that because of Jacobi's inversion theorem if $D$ is a divisor such that $r(-D)>0$ there is always a positive divisor $E$ of degree $g(X)-1$ and $D\equiv E.$

Now assume that $deg E=g(X)-1.$ Apply Riemann Roch and conclude that: $r(-D)=i(-D).$ Choose a base point $z_0$ on $X$ and let $u:X\mapsto Jac(X)$ the standard mapping from $X$ into its Jacobian. Let $K_{z_0}$ be the Riemann constant. Then Using Riemann vanishing theorem for theta functions we have the following non vanishing criteria for theta functions:
\begin{lem}
Let $D,$ $degD=g-1$ be a non positive divisor such that $r(-D)=0$ then $\theta\left(u(D)+K_{z_0}\right)\neq 0.$
\end{lem}

\section{Cyclic covers}
Let $\phi:X\mapsto \mathbb{CP}^1$ be a cyclic cover of the sphere of order $N$ prime number ramified above $m$ points $\lambda_1...\lambda_m$. Assume that $\lambda_i \ne \infty$ and let $P_i$ be the ramification point above $\lambda_i.$ Riemann Hurwitz formula applies $g(X)=\frac{(n-1)(m-2)}{2}.$ It is easy to see that $X$ satisfies the equation : $$y^N=\prod_{i=1}^m \left(x-\lambda_i\right)^{R_i}.$$  and $R_i \in \left\{1,2...N-1\right\}.$ $X$ has a cyclic group of automorphisms of order $N.$ An explicit generator of this group is: $T(y,x)=(\omega y,x),\omega^N=1.$ since $\lambda_i\neq \infty$  conclude: $\sum_{i=1}^m R_i=0 \bmod N.$ For each $j \in \left\{1...N-1\right\}$ define $t_j$ to be the number of $R_i$ such that $R_i=j.$ Then $\sum_{i=1}^{N-1}t_j=m. $ Associate to $X$ a vector $\alpha \in \mathbb{Z}^m,$ such that: $\alpha=\left(1...1,2...2,3...3...,N-1...N-1\right).$
The index $j$ appears $t_j$ times. Let $\alpha_i$ be the $i-th$ of $\alpha$. Since $\sum_{i=1}^m R_i=0\bmod N$ $\infty$ has precisely $N$ pre images in $X.$
Let these images be: $\left(\infty_1,...\infty_N\right)=\phi^{-1}\left(\infty\right).$
\begin{defn}
For $\xi$ is a divisor of degree $1$ on the sphere which is not a ramification point define $a\xi=\sum_{i=1}^N \phi^{-1}\xi.$
Extend the map to a divisor of any degree on the sphere
\end{defn}
On $X$ select a normalized homology basis $a_i,b_j$ and the set of the normalized canonical differentials $\omega_i.$  Choosing a base point $z_0 \in X$ define the mapping: $u:X\mapsto J(X).$
\begin{defn}
For the base point $z_0$ define the divisor:
${g_N}^1=\sum_{i=0}^{N-1} T^i\left(z_0\right).$
\end{defn}
\begin{defn}
Let $G_0 \in Jac(X)$ be a point such that $NG_0=u\left({g_N}^1\right).$
\end{defn}
Since $NP_i-{g_N}^1=0$ in the Jacobian conclude that: $u(P_k)=C_k+G_0$ and $NC_k=0.$ Let $\Delta$ be the canonical class and let $K_{z_0}$ be the Riemann constant i.e. $-2K_{z_0}=u(\Delta).$ Then $u(\Delta)=(N-1)\sum_{i=1}^m u(P_i)-2u\left({g_N}^1\right).$ Using the definition of $C_k$ rewrite the last expression as:$(N-1)\sum_{i=1}^m C_i+2(g(X)-1)G_0.$ Therefore Riemann's constant $K_{z_0}$ equals to:
$$
K_{z_0}=-\frac{N-1}{2}\sum_{i=1}^m C_i-(g(X)-1)G_0+E_2
$$
and $E_2$ is a point of order $2.$ i.e. $2E_2=0.$ Let $E_1=E_2+\frac{N-1}{2}\sum_{i=1}^mC_i$ then $E_1=-K_{z_0}-(g(X)-1)G_0.$
if $N$ is an odd number conclude that $2NE_1=0$ in J(X).
We like to formulate the main theorem which is the adaptation of [Ac2] p.26. This describes the vanishing order of theta functions at certain points of the Jacobian. Let $r$ be a total ramification of $f$ Note that: $r=m(N-1).$

\begin{thm} \label{one}
Let $\beta=\left(\beta_1...\beta_m\right) \in \left(0,...,N-1\right)$ such that $\sum \beta_j-\frac{r}{2}=0 \bmod N$ ($r$ is always even.)
Define a sequence of $N$ numbers $\tau_0...\tau_{N-1}$ satisfying the equations:
$$
  \sum_{i=1}^m\overline{\beta_i+k\alpha_i}=\frac{r}{2}-N\tau_k, 0\leq k\leq N-1
$$
Then the order of the theta function vanishing on the point $\sum_{i=1}^{m}\beta_iC_i-E_1$ in $Jac(X)$ is $\sum_{i=0}^{N-1} Max(0,\tau_i).$
\end{thm}
\begin{rmk}\label{rmk1}
Since $\overline{\beta}=\beta,$ notice:
$$
  \sum_{i=1}^m{\beta_i}=\frac{r}{2}-N\tau_0
$$

\end{rmk}
\noindent \textbf{Proof:\\}
The theta function vanishes on the point $\sum_{i=1}^{m}\beta_iC_i-E_1$ if and only if there is a positive divisor of degree $g-1$, $\psi$ such that:
$\sum_{i=1}^{m}\beta_iC_i-E_1=u(\psi)+K_{z_0}.$ Use the definition of $E_1$ and $C_i$ and the formula: $g(X)-1=r/2-N$ to write the last equality as:
$$
  u(\psi)=\sum_{i=1}^{m}\beta_iC_i-E_1-K_{z_0}=\sum_{j=1}^m\beta_ju(P_j)+(\tau_0-1){g_N}^1
$$
Where $\sum_{i=1}^{m}\beta_i=\frac{r}{2}-\tau_0N.$ by ~\ref{rmk1}. Let: $$D_1=\sum_{i=1}^{m}\beta_iP_i+(\tau_0-1){g_N}^1,$$
$D_1$ is a divisor (\textit{not necessarily positive}). Its degree given by the next proposition:

\begin{prop}
$D_1$ has degree $g(X)-1.$
\end{prop}

\noindent \textbf{Proof:\\}
$deg(D_1)=\sum_{i=1}^m\beta_i+N\times(\tau_0-1)$ By definition of $\tau_0$: $$\deg(D_1)=r/2-N\times \tau_0+N\times \tau_0-N$$ Recall that: $2g(X)-2=r-2N,$ hence $g-1=r/2-N$ and the proposition follows.$\blacksquare$\\\\
since $D_1$ is invariant under $T,$ $\forall f \in H^0\left(X,\mathcal{O}\left(-D_1\right)\right),$
$Tf \in H^0\left(X,\mathcal{O}\left(-D_1\right)\right).$
 $T$ is cyclic and unitary therefore $H^0\left(X,\mathcal{O}\left(-D_1\right)\right)$ has  a decomposition:   $H^0\left(X,\mathcal{O}\left(-D_1\right)\right)=\bigoplus L_\chi$ and $L_\chi$ is the vector space of $T$ eigenvectors with a character: $\chi:Z_n \mapsto \mathbb{C}.$ If  $N_\chi=dim L_\chi$ then Riemann Theorem applies that the order of vanishing of the theta function is: $N=\sum N_\chi.$ We attempt to find $N_\chi.$  $T$ is cyclic hence its characters are of the form $\omega^k$ for some $k,$ and  $\omega^N=1.$ Now $Ty=\omega y,$ and $Ty^k=\omega^ky^k.$ Conclude that if $f\in N_\chi$ then $f/y^k$ for some $k$ is a pull back of a function $g$ on $\mathbb{CP}^1.$ But $f/y^k\in H^0\left(X,\mathcal{O}\left(-D_1-div(y^k)\right)\right)$. Further $f/y^k$ corresponds to the functions that are pullbacks from the functions on the $\mathbb{CP}^1$ in the space: $H^0\left(X,\mathcal{O}\left(-D_1-div(y^k)\right)\right).$ Let ${V_k}^0$ be the space of $f\in H^0\left(X,\mathcal{O}\left(-D_1-div(y^k)\right)\right)$ that are pullbacks from functions on the sphere.
\begin{lem}
There is a divisor $\sigma_0$ with support on $\mathbb{CP}^1$ such that: $H^0\left(\mathbb{CP}^1,\mathcal{O}\left(-\sigma_0\right)\right)$ is isomorphic to ${V_k}^0.$
\end{lem}
\noindent \textbf{Proof}:\\
For a ramification point $\lambda_j$ let $\gamma_j=\left[\frac{kR_j+\beta_j}{N}\right]\times N$(i.e. $\gamma_j$ is the maximal number such that $\gamma_j\leq kR_j+\beta_j$ and $\gamma_j=0 \bmod N.$). Let $Q_{\lambda_i}$ be the point on $\mathbb{CP}^1$ that corresponds to $\lambda_i.$
Define $$\sigma_0=\sum_{j=1}^m\frac{\gamma_j}{N}Q_{\lambda_j}+(\tau_0-1)\phi(z_0)-k\frac{\sum_{j=1}^m R_j}{N}\infty.$$
We show that $\sigma_0$ is the desired divisor. Let $h:\mathbb{CP}^1\mapsto \mathbb{CP}^1$ be a function such that $divh\geq -\sigma_0.$ assume that $\hat{h}$ is a lift of $h$ to $X.$ Then $$div\hat{h}\geq-\gamma_jP_j-(\tau_0-1)g_1^N+k\frac{\sum_{j=1}^mR_j}{N}\sum_{i=1}^N\infty_i.$$ but $\gamma_j\leq kR_j+\beta_j,$ and from the definition of ${V_k}^0$ conclude, $\hat{h} \in {V_k}^0.$ Now assume that $\hat{f},$ a lift of a function $f$ on the sphere and $\hat{f} \in {V_k}^0.$ By definition at a point $P_j$ the $Ord_{P_j}(\hat{f})\geq -\beta_j-kR_j.$ $\hat{f}$ is a lift of $f$ hence,  $Ord_{P_j}(\hat{f})=0 \bmod N$ conclude that $Ord_{P_j}(\hat{f})\geq \gamma_j$ Or $Ord_{Q_{\lambda_j}}(f)\geq \frac{-\gamma_j}{N}.$ In other points of its support  The order of $div(-D_1-y^k)$ is divisible by $N.$ Thus: $f\in H^0\left(\mathbb{CP}^1,\mathcal{O}\left(-\sigma_0\right)\right).$
$\blacksquare$\\\\
The immediate conclusion from the lemma is that: $dim{V_k}^0=dim H^0\left(\mathbb{CP}^1,\mathcal{O}\left(-\sigma_0\right)\right).$ Let us compute the degree of $\sigma_0.$ By definition of $\sigma_0$ we have:
$$deg\sigma_0=\frac{1}{N}\times\left(\sum_{i=1}^m\gamma_i+(\tau_0-1)N-\sum_{j=1}^mkR_i    \right)$$
But $\gamma_i=\beta_i+kR_i-\overline{\beta_i+kR_i}$ by the definition of $\gamma_i.$ Substituting this expression into $\gamma_i$ rewrite the last expression as: $$\frac{1}{N}\times\left(\sum_{i=1}^m\beta_i+kR_i-\left(\overline{\beta_i+kR_i}\right)+(\tau_0-1)N-\sum_{j=1}^mkR_i    \right)$$ Cancel $kR_i$ and apply the definition of $\tau_0$ to simplify further: $$\frac{1}{N}\times\left(\frac{r}{2}-\sum_{i=1}^m\left(\overline{\beta_i+kR_i}\right)\right)-1.$$ By definition of $\tau_i$ this equals to: $\tau_k-1.$ Apply  Riemann Roch on $\mathbb{CP}^1$ to conclude that $dimH^0\left(\mathbb{CP}^1,\mathcal{O}\left(-\sigma_0\right)\right)=Max(\tau_i,0).$
$\blacksquare$\\
\noindent The discussion in section 2 produces the following corollary:
\begin{cor}
Choose $\beta_i$ as in theorem but $\tau_i=0$ then under the conditions of last theorem $\theta\left(\beta_iC_i-E_1\right)$ is not vanishing.
\end{cor}
Recall that $-E_1-K_{z_0}=-(g(X)-1)G_0,$ Or $-E_1=K_{z_0}+(g(X)-1)G_0.$ Rewrite the divisor from corollary as:
$$
\sum_{i=1}^m \beta_i(P_i-G_0)+K_{z_0}+(g(X)-1)G_0=\sum_{i=1}^m \beta_iP_i-\beta_iG_0+K_{z_0}+(g(X)-1)G_0=$$
$$\sum_{i=1}^m\beta_iP_i-\frac{r}{2}G_0+(\frac{r}{2}-N)G_0 +K_{z_0}$$
and the last expression is readily seen to be equal to :
$\sum_{i=1}^m\beta_iP_i+K_{z_0}-\sum_{i=1}^N \infty_i.$
\begin{cor} \label{cornonvanishing}
Let $\beta_i$ be selected such that $\tau_i=0, 0\leq i\leq N-1 $ then $$\theta\left[u(\sum_{i=1}^m \beta_iP_i)+ K_{z_0}-u(\sum_{j=1}^N\infty_j)\right](0,\tau)\neq 0.$$
\end{cor}
\begin{rmk}
Gabino in [GG] obtained similar results but the theorem stated here seems to be stronger and was independently obtained. See also [EF] for an alternative proof where non positive divisors of degree $g-1$ are replaced with the more traditional special divisors of degree $g.$
\end{rmk}

\section{N=3 example}
We work out the general $N=3$ example following [Ac1]. Let us represent the curve as : $y^3=\prod_{i=1}^s(x-p_i)\prod_{j=1}^t(x-q_j)^2.$ and $s+2t=0\bmod 3.$ Then $G_0$ be a point satisfying $3G_0=u(z_1+z_2+z_3)$ in the Jacobian. The ramification points lying above $p_i$ and $q_j$ will be respectively: $a_k=A_k+G_0$ and $b_k=B_k+G_0.$ Consider the sum $\sum_{i=1}^s \epsilon_i A_i+\sum_{j=1}^t\delta_jB_j-E_1.$ In the vector $\alpha=(1...1,2...2)$ the ones appearing $s$ times and $2$ appearing $t$ times.
Then we have the following conditions on $\epsilon_i,\delta_j:$
\begin{itemize}
  \item $\sum\epsilon_i + \sum \delta_j = s+t-3\tau_0$
  \item $\overline{\sum(\epsilon_i+1)} + \overline{\sum (\delta_j+2)} = s+t-3\tau_1$
  \item $\overline{\sum(\epsilon_i+2)} + \sum (\overline{\delta_j+1)} = s+t-3\tau_2$
\end{itemize}
Now rewrite the period as: $E_1-\sum_{S_1}A-2\sum_{S_2}A-\sum_{T_1}B-\sum_{T_2}2B$ where $S_1,S_2$ are subsets where appearing $1$ and $2$ in the $A$ part of the sum and $T_1$ and $T_2$ appearing in the $B$ part of the sum. Accordingly $|S_i|=s_i$ and similarly $|T_i|=t_i.$ Finally $S_0,T_0$ be subsets of indices such that $\epsilon_i=\delta_i=0.$ Define $\mu_0=s_0-t_2$, $\mu_1=s_1-t_1$ and $\mu_2=s_2-t_0.$  Write the condition on $\tau_i$ as:
$3\tau_0=\mu_0-\mu_2, 3\tau_1=\mu_2-\mu_1, 3\tau_2=\mu_1-\mu_0.$ Therefore  $\mu_0=\mu_1=\mu_2=\frac{t-s}{3}$ guarantees non vanishing. We showed:
\begin{thm}
Let $\left\{S_0,S_1,S_2\right\},\left\{T_0,T_1,T_2\right\}$ be a partition of $a_i,1\leq\ i \leq s $ and $b_j,1\leq j\leq t$ respectively. if $t_2=s_0+(t-s)/3,
t_1=s_1+(t-s)/3, \tau_0=s_2+(t-s)/3.$ Then $\theta$ does not vanish on the following point of the Jacobian: \\$u\left((\sum_{S_1}A+2\sum_{S_2}B+\sum_{T_1}B+2\sum_{T_2}B)-\infty_1-\infty_2-\infty_3\right)+K_{z_0}.$
\end{thm}
\section {The non singular case}
As a second example of applying the result assume $R_i=1.$ Then $m=pN.$
Then we can identify the coefficients $\beta_i$ with a vector $v$ in the integral lattice: $Z^{m}$ such that the $i$-th coordinate of $v$ is $\beta_i.$ Now the vector $\alpha=\left(1,1,1...1,1)\right)$ and consequently the $\tau_k$ are defined as:
$$\sum_{i=1}^{m}\overline{\beta_i+k}=\frac{r}{2}-N\tau_k.$$ since $r=m\left(N-1\right)$ we can rewrite the last expression as :
$$
\sum_{i=1}^m\overline{\beta_i+k}=m\frac{N-1}{2}-N\tau_k
$$
Let $t_l$ be the number of times that $\beta_i=l,0<l<N-1.$ w.l.o.g we can assume that $max\left(t_1...t_{N-1}\right)=t_1$ then:
$$
\sum_{i=1}^{m}\overline{\beta_i+k}=N\frac{N-1}{2}t_1+2(t_2-t_1)+3(t_3-t_1)+...(N-1)(t_{N-1}-t_1)
$$
Because $t_i-t_1\leq 0$ The condition for non vanishing turns out to be:
$t_i=t_1 \forall i$ and $t_1=p.$
We obtained the following theorem:
\begin{thm} \label{two}
Assume that the $R_i=1, \forall i$ Then the divisor $\sum_{i=1}^N\beta_iP_i+K_{z_0}-\sum_{i=1}^N \infty_i$ is non vanishing if and only if the number of $\beta_i$ such that $\beta_i=l, 0\leq l\leq N-1$ is $k.$
\end{thm}
\section{Properties of divisors}
We examine the properties of divisors appearing in \textbf{Corollary} \textbf{~\ref{cornonvanishing}}. This will be useful in the sequel:
\begin{prop}
Let $D=\sum_{i=1}^m \beta_iP_i$ and let $E=\sum_{i=1}^m\left(\overline{\beta_i+k\alpha_i}\right)P_i.$ Then $D\equiv E$
\end{prop}
\noindent \textbf{Proof:\\} We show that there exists a function on $X$ such that its divisor is: $D/E.$ Define: $$f=\frac{\prod_{i=1}^m (x-\lambda_i)^{\beta_i/N+k\alpha_i/N}}{\prod_{i=1}^m(x-\lambda_i)^{\overline{\beta_i+k\alpha_i}/N+k\alpha_i/N}}$$
Clearly the divisor of $f$ is precisely $D/E.$ We verify that this is indeed a function in $X.$ Using the definition of $y$ we rewrite $f$ as :$$y^{-k}\times \frac{\prod_{i=1}^m (x-\lambda_i)^{\beta_i/N+k\alpha_i/N}}{\prod_{i=1}^m(x-\lambda_i)^{\overline{\beta_i+k\alpha_i}/N}}   $$ and by definition $\beta_i+k\alpha_i-\overline{\beta_i+k\alpha_i}=h_iN,  h_i\in \mathbb{Z}. $ Hence $f=y^{-k}\times \prod_{i=1}^m\left(x-\lambda_i\right)^{h_i}$ $\blacksquare$\\
\begin{prop}
Let $D = u(\sum_{i=1}^m\beta_iP_i)+K_{z_0}-u(\sum_{i=1}^N\infty_i)$ as in ~\ref{cornonvanishing}. Then $-D=u(\sum_{i=1}^m(N-1-\beta_i))+K_{z_0}-u(\sum_{i=1}^N \infty_i)$
\end{prop}
\noindent \textbf{Proof:\\} We can write:
$$
 -D=-u(\sum_{i=1}^m\beta_iP_i)-K_{z_0}+u(\sum_{i=1}^N \infty_i)=
$$
$$
-u(\sum_{i=1}^m\beta_iP_i)-2K_{z_0}+ K_{z_0}+2\sum_{i=1}^N \infty_i-\sum_{j=1}^N\infty_i
$$
Because of the definition of $K_{z_0}$ we obtain that $-2K_{z_0}\equiv u\left(\sum_{i=1}^m(N-1)P_i-2\sum_{j=1}^m{\infty_i}\right).$ Hence:
$$-\sum_{i=1}^m\beta_iP_i-2K_{z_0}+K_{z_0}+2u(\sum_{i=1}^N \infty_i)-u(\sum_{i=1}^N\infty_i)=$$$$
\sum_{i=1}^m(N-1-\beta_i)P_i+K_{z_0}-\sum_{j=1}^N \infty_i
$$
$\blacksquare$\\
\begin{lem}
For $0\leq k\leq N-1$ Let  $$f_k(z) = \prod_{i=1}^m\left(z-\lambda_i\right)^{\frac{\overline{\beta_i+k\alpha_i}-(\frac{N-1}{2})}{N}}\sqrt{dz}
$$
Then $f_k(z)$ is a meromorphic whose divisor is equivalent to: $\sum_{i=1}^m\beta_iP_i-\sum_{i=1}^N \infty_i$
\end{lem}
\noindent \textbf{Proof:\\} The order of $z-\lambda_i$ at $P_i$ is  $N.$ Hence the order of the of $(z-\lambda_i)^{\frac{\overline{\beta_i+k\alpha_i}-(\frac{N-1}{2})}{N}}\sqrt{dz}$ is exactly $\overline{\beta_i+k\alpha_i}.$ at $\infty_i$ the order is: $\sum_{i=1}^m\overline{\beta_i+k\alpha_i}-\frac{r}{2}-1=-1$ and hence the divisor is: $$\sum_{j=1}^m \overline{\beta_j+k\alpha_j}P_j-\sum_{i=1}^N \infty_i$$ which is equivalent to : $\sum_{j=1}^m\beta_jP_j-\sum_{i=1}^N \infty_i.$
\noindent $\blacksquare$
\section{Algebraic construction of the Szego Kernel}
Let us recall the definition of the Szego Kernel.
\begin{defn}
For $e\in \mathbb{C}^g.$ if $\theta[e]\neq 0$ define the Szego kernel by the following equation: $$S(P,Q|e)=\frac{\theta[e]\left(u(P-Q)\right)}{\theta[e](0,\tau)E(P,Q)},P,Q\in C.$$
\end{defn}
$E(P,Q)$ is the prime form. $e$ depends only on its in the Jacobian , $J(X)$. $R(P,Q|e)$ has the following properties that are well known [F] (p.19,p.123),[EG2] (Proof of Theorem 4.7):

\begin{itemize}
  \item $S[e](P,Q)$ is a $\left(\frac{1}{2},\frac{1}{2}\right)$ form with a simple pole along the diagonal
  \item $S[e]\left(P,Q\right)$ has divisor $[e-K_{z_0}]$ with respect to variable $Q$
  \item $S[e]\left(P,Q\right)$ has a divisor $[-e-K_{z_0}]$ with respect to variable $P.$
  \item $S[e](P,Q$ is a unique up to a constant $\left(\frac{1}{2},\frac{1}{2}\right)$ form that satisfies the previous properties
\end{itemize}

We generalize the approach of Nakayashiki to give the following expression to the Szego kernel:
\begin{thm} \label{szego}
Let $P=(x_1,y_1),Q=(x_2,y_2)\in C.$ Choose
 $\overline{\beta}=\left(\beta_1...\beta_m\right) \in \mathbb{Z}^m$ be as in \textbf{~\ref{cornonvanishing}}. if  $e = \sum_{i=1}^mu(\beta_iP_i)+K_{z_0}-u(\sum_{i=1}^N\infty_i)$ Let \begin{equation}\label{szegoeq}
 F_{\overline{\beta}}\left(P,Q\right) =\frac{1}{N}\frac{\sum_{k=0}^{N-1}\prod_{i=1}^m\left(x_1-\lambda_i\right)^{\frac{\overline{\beta_i+k\alpha_i}-\frac{N-1}{2}}{N}}\times \prod_{i=1}^m\left(x_2-\lambda_i\right)^{\frac{\overline{-\beta_i-k\alpha_i}+\frac{N-1}{2}}{N}}\sqrt{dx_1}\sqrt{dx_2}}{x_2-x_1}
\end{equation}
Then \begin{equation}
S[e](P,Q)=F_{\widetilde{\beta}}\left(P,Q\right)
\end{equation}

\end{thm}
\noindent \textbf{Proof:\\} We verify that $F_{\widetilde{\beta}}\left(P,Q\right)$ satisfies the properties characterizing $S[e]\left(P,Q\right).$ the RHS of the equation  \textbf{(~\ref{szegoeq})}.

   $F_{\widetilde{\beta}}\left(P,Q\right)$ is regular outside $P=Q.$ If $P\neq Q$ we need to check the case when $x_1=x_2$ and $y_1\neq y_2.$ This means that $y_2=\omega^jy_1.$ and $\omega^N=1.$\\ Now, $\overline{\beta_i+k\alpha_i}=\beta_i+k\alpha_i-h_{ki}N.$ Rewrite $F_{\widetilde{\beta}}\left(P,Q\right)$ as:
$$\frac{1}{N}\times\frac{1}{x_2-x_1}\times \sum_{k=0}^{N-1}\left(\frac{y_1}{y_2}\right)^k\times\left(\frac{x_1-\lambda_i}{x_2-\lambda_i}\right)^{h_{ik}}\prod_{i=1}^m{\left(\frac{x_1-\lambda_i}{x_2-\lambda_i}\right)}^{\frac{\beta_i-\frac{N-1}{2}}{N}}\sqrt{dx_1}\sqrt{dx_2}$$
if $y_2=\omega^k y_1.$ In the limit when $x_1\rightarrow x_2$,  $y_2\rightarrow\omega^k y_1.$ and  $\frac{\sum_{k=0}^{N-1}\left(\frac{y_1}{y_2}\right)^k}{x_2-x_1}\rightarrow 0.$ Therefore $F_{\widetilde{\beta}}(P,Q)$ is regular when $x_1=x_2$ but $y_1\neq y_2.$

\noindent
Let us calculate the expansion of $$\frac{1}{N}\sum_{k=0}^{N-1}\prod_{i=1}^m\left(x_1-\lambda_i\right)^{\frac{\overline{\beta_i+k\alpha_i}-\frac{N-1}{2}}{N}} \times \prod_{i=1}^m\left(x_2-\lambda_i\right)^{\frac{\overline{-\beta_i-k\alpha_i}+\frac{N-1}{2}}{N}}$$ as a function of $x_2$ when the expansion is around $x_1.$ Assuming $x_1=x_2$ we get that the leading coefficient is $\frac{1}{N}\sum_{i=0}^{N-1} 1=1$ Now for the coefficient in $x_2-x_1$ we obtain using the derivative product rule that the coefficient is:
\begin{equation}
\frac{1}{N}\sum_{i=0}^{N-1}\sum_{j=1}^m\left(\overline{-\beta_j-k\alpha_j}/N+\frac{(N-1)}{2N}\right)\times\frac{1}{x_2-\lambda_j}=
\frac{1}{N}\sum_{j=1}^{m}\sum_{i=0}^{N-1}\left(\overline{-\beta_j-k\alpha_j}/N)+\frac{r}{2}\right)\times\frac{1}{x_2-\lambda_j}
\end{equation}
But $\sum_{j=1}^{m}\left(\overline{-\beta_j-k\alpha_j}/N\right)+\frac{r}{2}=0$ and hence
$$\frac{1}{N}\sum_{j=1}^{m}\sum_{i=0}^{N-1}\left(\overline{-\beta_j-k\alpha_j}/N)+\frac{r}{2}\right)\times\frac{1}{x_2-\lambda_j}=0$$ as well.
Taking the second derivative according to $x_1.$ to calculate the coefficient of $(x_1-x_2)$ we arrive to the following result:
\begin{prop}
The expansion of  $F_{\widetilde{\beta}}\left(P,Q\right)$ around $P$ a non branch point is:
\begin{equation}
F_{\widetilde{\beta}}\left(P,Q\right) = \frac{\sqrt{dx_1}\sqrt{dx_2}}{x_2-x_1}\left[1+\frac{1}{2N}\sum_{i,j=1}^{i,j=m}\frac{q(\beta_i,\beta_j)}{(x_2-\lambda_i)(x_2-\lambda_j)}\times(x_1-x_2)^2+...\right]
\end{equation}
where $$q(\beta_i,\beta_j)=\sum_{k=0}^{N-1}\left(\left\{(\beta_i+kR_i)/N\right\}-\frac{N-1}{2N}\right)\times\left(\left\{(\beta_j+kR_j)/N\right\}-\frac{N-1}{2N}\right)$$
Where $x_1,x_2$ are the local coordinate around $P,Q$ respectively.
\end{prop}
To complete the proof of theorem \textbf(~\ref{szego}) note that $L(P,Q)=F_{\overline{\beta}}(P,Q)-S[e](P,Q)$ are a section of a line bundle $L_{[e-K_{z_0}]}\bigotimes L_{[-e-K_{z_0}]}.$ Because of the expansion of $F_{\overline{\beta}}(P,Q)$ conclude that $L(P,Q)$ is a holomorphic section of the line bundle. But $$H^0\left(L_{[e-K_{z_0}}]\bigotimes L_{[-e-K_{z_0}]}\right)=  H^0\left(L_{[e-K_{z_0}]}\right)\bigotimes H^0\left(L_{[-e-K_{z_0}]}\right)=0$$
and thus $F_{\overline{\beta}}(P,Q)=S[e](P,Q)$ as required.
$\blacksquare$
\begin{rmk}
The above argument is exactly the method adopted in [Na] to prove the claim for the non singular case. See [EG] for a slightly different approach.
\end{rmk}
Based on the the formula given at the beginning of the section [Na] shows the following expansion for $S[e](P,Q)$ in terms of theta functions:
\begin{cor}
The expansion of the Szego kernel can be given in terms of theta functions as follows:
$$S[e](P,Q)=\frac{\sqrt{dx_1}\sqrt{dx_2}}{x_1-x_2}\times\left[1+\sum_{i=1}^g\frac{\partial \log\theta[e]}{\partial z_i}(0)u_i(x_1)(x_1-x_2)+...\right]$$
$u_i(x)$ is the coefficient of $dx_1$ in the expansion of the holomorphic $v_i(x).$
\end{cor}
Comparing the expansions conclude the following result:
\begin{cor}\label{thetadervanish}
$$\frac{\partial \theta[e]}{\partial z_i}(0)=0$$
\end{cor}
The following is obtained by multiplying the expansions:
\begin{lem}
\begin{equation}
S[e](P,Q)S[-e](P,Q)=\frac{dx_1dx_2}{(x_1-x_2)^2}\left[1+\frac{1}{N}\sum_{i,j=1}^m \frac{q\left(\beta_i,\beta_j\right)}{(x_2-\lambda_i)(x_2-\lambda_i)}(x_1-x_2)^2+...\right]
\end{equation}
\end{lem}
\section{Algebraic construction for the canonical differential}
We construct the canonical differential algebraically for cyclic covers.
\begin{defn}
The canonical symmetric differential is a $\omega(x,y)$ is a meromorphic  one differential with respect to $x,y \in C$, having a unique pole of second order when $z$ tends to $w$ with a leading expansion coefficient of $1.$ Further for a canonical homology basis $a_i,b_j, 1\leq i,j\leq g$ we have:
$$
\int_{a_i}\omega(x,y)=0
$$
for fixed $y.$
\end{defn}
First we remind the reader of a possible basis for the holomorphic differentials on
$C.$
\begin{lem}
Let $s_l(z)=\prod_{i=1}^m(x-\lambda_i)^\frac{\overline{lR_i}}{N}., 0\leq l\leq N-1$
Then a basis for holomorphic differentials is given by: $$\frac{z^{j-1}dz}{s_l(z)},$$ where $j=1...d(l), d(l)=Max\left(\sum_{i=1}^m\frac{\overline{lR_i}}{N}-1,0\right).$
\end{lem}
\noindent \textbf{Proof:\\}
The order of $z-\lambda_i$ at $\lambda_i$ is $N.$ Hence the order of $(z-\lambda_i)^{\frac{\overline{lR_i}}{N}}$ is $\overline{lR_i}<N-1.$ Hence we have non trivial $0$ at $\lambda_i.$ For $\infty_i, i=1...n.$ The order of $y_l(z)$ is $\sum_{i=1}^N\frac{\overline{lR_i}}{N}.$ Thus if $j<\sum_{i=1}^N\frac{\overline{lR_i}}{N}$ we will not have a pole at $\infty_i.$
$\blacksquare.$\\
Let,  $$P_l^{(l)}\left(z,w\right)=\sum_{n=0}^{d(l)+1}A_n^{(l)}(w)(z-w)^n,$$ such that:
\begin{enumerate}
\item $A_0^{(l)}\left(w\right)=\prod_{i=1}^m\left(w-\lambda_i\right)$
\item $A_1^{(l)}=\sum_{i=1}^m\overline{lR_i}/N\times\frac{A_0^{(l)}\left(w\right)}{w-\lambda_i}$
\item $deg_{w}P_l^{(l)}\leq d(N-l)+1$
\end{enumerate}
Set:
\begin{equation}\xi_0(x,y)=\frac{dz(x)dz(y)}{\left(z(x)-z(y)\right)^2}
\end{equation}
\begin{equation}\xi_l\left(x,y\right)=\frac{P_l^{(l)}\left(z(x),z(y)\right)dz(x)dz(y)}{s_l(x)s_{N-l}(y)\left(z(x)-z(y)\right)^2}
\end{equation}
\begin{equation} \xi\left(x,y\right)=\frac{1}{N}\sum_{i=0}^{N-1}\xi_l\left(x,y\right)
\end{equation}
\begin{prop}
\begin{enumerate}
\item $\xi\left(x,y\right)$ is holomorphic outside the diagonal set $\left\{x=y\right\}.$
\item For a non branch point $P\in X$ take $z$ to be a local coordinate around $P.$ Then the expansion in $z(x)$ at $z(y)$ is :
    $$
    \xi\left(x,y\right)=\frac{dz(x)dz(y)}{\left(z(x)-z(y)\right)^2}+
     O \left(\left(z(x)-z(y)\right)^0\right)
    $$
\end{enumerate}
\end{prop}
\noindent \textbf{Proof:\\}
To show the proposition we need first need to show that if $z(P)=z(Q)$ but $P\neq Q$ on $X,$ then $\xi(P,Q)$ is still non singular. Assume that $P=(p_1,q_1)$ and $Q=(p_1,\omega^rq_1).$ Let us examine the leading term of the expansion of $\xi_l(x,y)$ around $Q.$ By definition of $\xi_l(x,y)$ it is:
$$\frac{A_0^{(l)}\left(z(Q)\right)}{\omega^{rl}{q_1^l/\prod_{i=1}^m(z-\lambda_i)^{[lR_i/N]}}q_1^{N-l}\times(z(P)-z(Q))^2}=$$
$$=\omega^{-rl}\times(z(P)-z(Q))^{-2},$$ ( by definition of $A_0^{(l)}$)
Then if $r=0$ (i.e. $P\rightarrow Q$) the leading coefficient of $\xi(x,y)$ is: $\frac{1}{\left(z(P)-z(Q)\right)^2}.$ if $r>0,$ summing we obtain that the coefficient is $0.$ Next we compute the coefficient of $ \frac{1}{z(P)-z(Q)}$in the expansion of $\xi_l(x,y),l>0.$ Apply the product rule for derivatives to obtain that the coefficient of $ \frac{1}{z(P)-z(Q)}$ is expanding around $Q$ is:
$$
\frac{\sum_{i=1}^m\left\{lR_i/N\right\}A_0^{\left(l\right)}(z(Q))}{z(Q)-\lambda_i}s_l(Q)-\sum_{i=1}^m\left\{lR_i/N\right\}\frac{s_l\left(z(Q)\right)}{z(Q)-\lambda_i}A_0^{(l)}(z(Q))=0
$$
therefore the coefficient of $\frac{1}{z(P)-z(Q)}$ is $0$ and the proposition is proved.\\
$\blacksquare$

As an immediate corollary of the proposition we have that:
\begin{cor}
$$\omega(x,y)-\xi(x,y)$$ is holomorphic on $X\times X.$
\end{cor}
Thus by the corollary there exist polynomials $P_k^{(l)}$ such that:
$$
\omega(x,y)-\xi(x,y)=\sum_{l=1}^{N-1}\sum_{k=1, k\neq l }^{N-1}\frac{P_k^{(l)}\left(z(x),z(y)\right)dz(x)dz(y)}{s_k(z(x))s_{N-l}(z(y))}
$$
Where by modifying the definition of $P_{l}^{(l)}$ we can exclude the terms $k=l.$ as before we can write $$P_k^{(l)}(z,w)=\sum_{j=0}^{d(k)}A_{kj}^{(l)}(w)(z-w)^j.$$
Note that $deg_wP_k^{(l)}(z,w)\leq d(N-l).$
Our aim is to show the following proposition:
\begin{prop}\label{temp}
\begin{equation}
\sum_{l=1}^{N-1}A_{l2}^{(l)}\left(\lambda_i\right)=-\prod_{j=1,i\neq j}^m(\lambda_i-\lambda_j)\frac{\partial}{\partial \lambda_i}\log\det C
\end{equation}
and $C$ is a $g(X)\times g(X)$ period matrix of non normalized form : $$\left(\int_{a_i}z^{j-1}dz/s_l(z)\right).$$
\end{prop}
\noindent \textbf{Proof:\\} Let us take a local $t=\left(z-\lambda_i\right)^{\frac{1}{N}}$ coordinate around $Q_i.$ Then we have that the condition that $\int_{a_j}\omega(x,y)=0$ is equivalent to the coefficient of the expansion around $Q_i$ in $dt,tdt,...t^{N-2}dt$ is vanishing. A short calculation shows that this is equivalent to
\begin{equation} \omega^{(l)}\left(x\right)=\frac{1}{N}\frac{P_l^{(l)}(z(x),\lambda_i)dz(x)}{s_l(x)(z(x)-\lambda_i)^2}+\sum_{k=1,k\neq l}^{N-1}\frac{P_k^{(l)}\left(z(x),\lambda_i\right)dz(x)}{s_k(x)}
\end{equation}
vanishing when we integrate around $a_j.$ Let us write this explicitly: Note,  $$\frac{\partial}{\partial \lambda_i}\frac{dz}{s_l}=\left\{lR_i/N\right\}\frac{dz}{s_l\left(z-\lambda_i\right)}$$ and
$$P_l^{(l)}(z,\lambda_i)=\left\{lR_i/N\right\}\prod_{j=1}^m\left(\lambda_i-\lambda_j\right)(z-\lambda_i)+\sum_{j=0}^{d(l)-2}A_{l,j+2}^{(l)}(z-\lambda_i)^{j+2}.$$
hence:
\begin{equation}
\frac{\prod_{j=1}^m\left(\lambda_i-\lambda_j\right)}{N}\frac{\partial}{\partial \lambda_i}\int_{a_h}\frac{dz}{s_l}+\frac{1}{N}\sum_{j=0}^{d(l)-1}{A^{(l)}_{l,j+2}(\lambda_i)}\int_{a_h}(z-\lambda_i)^jdz/s_l
\end{equation}
$$+\sum_{k=1,k\neq l}\sum_{j=0}^{d(k)-1}{A^{(l)}_{k,j}}(\lambda_i)\int_{a_h}(z-\lambda_i)^jdz/s_k=0
$$
Following [BR],[Na](see also [EG] for a slightly different approach.) For a fixed $l$ regard the equations above as $g(X)$ equations in $g(X)$ variables, $A^{(l)}_{k,r}.$ The matrix of these equations is the $g(X)\times g(X)$ matrix $B,$ $$B=\int_{a_h}\left(z-\lambda_i\right)^{j-1}dz/s_l(z), l=1...N-1,j=1...d(l).$$ For each $l$ define matrices $B_l$ obtain from $B$ by replacing the column $\int_{a_h}\frac{dz}{s_l(z)}, 1 \leq g(X)$ with the column:$\frac{\partial}{\partial \lambda_i}\int_{a_h}\frac{dz}{s_l(z)}.$ Then
by Cramer's rule: $$A_{l2}^{(l)}=\frac{\det B_l}{\det B}.$$ Expand $(z-\lambda_i)^k, k\geq 0$ using the binomial formula and perform elementary operations on columns to obtain that $\det B=\det C.$ We now prove a similar proposition for the matrix $\det B_l.$
\begin{prop}
$$\det B_l=\sum_{i=0}^{d(l)}\det C_i$$ where $C_i$ is the matrix $C$ where the $i-th$ column is replaced with $\frac{\partial}{\partial \lambda_i}\frac{z^idz}{s_l(z)}.$
\end{prop}
\textbf{Proof:}
Assume inductively that the proposition is true for any minor of $B_l$ that contains the first $v$ columns.
The polynomials $1,(z-\lambda_i),...(z-\lambda_i)^{h-1}$ are basis for the polynomials $p(z), deg_{z}p(z)\leq h-1.$  write: $$(z-\lambda_i)^h=(z-\lambda_i)^h-z^h+z^h-{\lambda_i^h}+\lambda_i^h
=\sum_{i=0}^{v=h-1}c_v{\left(z-\lambda_i\right)}^v+z^h-\lambda_i^h
$$
Thus the $\det B_l$ is equal to the determinant of a matrix where the last column is replaced by $\int_{a_i}\left(z^{d(l)}-\lambda_i^{d(l)}\right)\times\frac{dz}{s_l}.$ Hence $\det B_l=\det D_l-\det E_l$ where the last columns of $D_l,E_l$ are replaced by: $z^{d(l)}$ and $\lambda_i^{d(l)}$ respectively. By induction assumption we see that $\det D_l=\sum_{v=1}^{d(l)-1}\det C_v.$ So to finish the proof we need to show that $\det E_l=-\det C_{d_l}.$ On the matrix $E_l$ perform the elementary operations by swapping the first and last column. dividing the first column by $\lambda_i^h$ and multiplying the last column by $\lambda_i^h$ respectively results in a matrix $E_l^{'}$ where the first $d(l)-1$ columns are $\int_{a_i}(z-\lambda_i)^{v}dz/s_l(z)$ and the last column is $\lambda_i^{d(l)}\frac{\partial}{\partial \lambda_i}\int_{a_h}\frac{dz}{s_l(z)}.$
Perform elementary operations on the first $d(l)-1$ to replace the columns with $\int_{a_i}\frac{z^vdz}{s_l}.$ Now we can write :

$$\frac{z^j}{z-\lambda_i}=\sum_{h=1}^{j-1}\lambda_i^{h-1}z^{j-h}+\frac{\lambda_i^j}{z-\lambda_i}.$$
But: $$\frac{\partial}{\partial \lambda_i}\frac{z^jdz}{s_l}=\left\{lR_i/N\right\}z^j\frac{dz}{s_l\left(z-\lambda_i\right)}$$
Hence:
$$
\frac{\partial}{\partial \lambda_i}\int_{a_k}\frac{z^jdz}{s_l}=\lambda_i^j\times \frac{\partial}{\partial \lambda_i}\int_{a_k}\frac{dz}{s_l}+\left\{lR_i/N\right\}\sum_{h=1}^{j-1}\lambda_i^{h-1}\int_{a_k}z^{j-h}\frac{dz}{s_l}
$$
Set $j=d(l)$ and using the fact ( yet again!) that elementary operations don't alter the determinant to conclude the result.
To finish the proof of \textbf{~\ref{temp}} observe that $$\frac{\sum_{h=1}^{N-1}\sum_{j=0}^{d(h)}\det C_i}{\det C}=\frac{\partial}{\partial \lambda_i}\log(\det C).$$\\
$\blacksquare$\\
Let us define the following object we will work closely when showing Thomae:
\begin{defn}
Let $P=(x,y) \in X$ be a non branch point with a local coordinate $z.$ Define: $$G_z(z)=\lim_{y \rightarrow x}\left[\omega(x,y)-\frac{dz(x)dz(y)}{{(z(y)-z(x))}^2}\right]
$$
\end{defn}
Now taking the local coordinate $t={(z-\lambda_i)}^{\frac{1}{N}}$ around the branch point $\lambda_i$ we have the following corollary:
\begin{cor}\label{atomic}
The coefficient of $t^{N-2}dt^2$ in the expansion of $G_z(z)$

\noindent in $t={\left(z-\lambda_i\right)}^{\frac{1}{N}}$ is:
$$-N\sum_{j=1,j\neq i}^m\frac{\gamma_{ij}}{\lambda_i-\lambda_j}-N\log \det C,$$ where
$$\gamma_{ij}=\sum_{h=0}^{N-1}\left\{hR_i/N\right\}\left\{hR_j/N\right\}$$
\end{cor}
To proceed further we learned the following from [Na] see([F] Corollary 2.12) that if $e$ belongs to the Jacobian such that $\theta[e](0,\tau)\neq 0$ then: $$
S[e]\left(x,y\right)=\omega(x,y)+\sum_{i,j=1}^g\frac{\partial^2\log\theta[e](0)}{\partial z_i\partial z_j}v_i(x)v_j(y),
$$
$v_i(x),v_j(x)$ are holomorphic differentials on the surface. Using this we obtain the following expression for $G_z(z)$
\begin{prop}
\begin{equation}
G_z(z)=\frac{1}{N}\sum_{i,j=1}^m\frac{q(\beta_i,\beta_j)dz(x)^2}{(z(x)-\lambda_i)(z(x)-\lambda_j)}-\sum_{i,j=1}^g\frac{\partial^2\log\theta[e](0)}{\partial z_i\partial z_j}v_i(x)v_j(y),
\end{equation}
\end{prop}
Passing to the local coordinate $t=\left(z-\lambda_i\right)^{\frac{1}{N}}$ we obtain the following corollary:
\begin{cor}\label{CorYaacov}
The coefficient of $t^{N-2}dt^2$ in the Laurent expansion of $G_z(z)$ is :
\begin{equation}\label{Formulamain}
2N\sum_{j=1,j\neq i}^{j=m}\frac{q(\beta_i,\beta_j)}{\lambda_i-\lambda_j}-\frac{1}{\left(N-2\right)!}\sum_{r,s=1}^g\sum_{\alpha=0}^{N-2}\frac{(N-2)!}{\alpha!(N-\alpha)!}\frac{\partial^2\log\theta[e](0)}{\partial z_i\partial z_j}{v_r}^{\alpha}(P_i){v_s}^{N-2-\alpha}(P_i),
\end{equation}
${v_r}^{\alpha}(P_i)$ is the coefficient of $dt$ in the expansion of $v_r(x)$ in the local coordinate $t.$
\end{cor}
\section{Variational formula for the period matrix and Thomae for general cyclic covers}
[Na] shows the following formula that can be generalized to any cyclic cover: \begin{thm}\label{variational}
If $\tau$ is a period matrix with respect to the fixed homology basis ${a_i,b_j},1\leq i,j\leq g$ then
\begin{equation}
\frac{d\tau_{jk}(0)}{dt}=\frac{1}{N(N-2)!}\sum_{\alpha=0}^{N-2}\left(
                                                                 \begin{array}{c}
                                                                   N-2 \\
                                                                   \alpha \\
                                                                 \end{array}
                                                               \right)
{v_j}^{\alpha}\left(Q_i\right)v_k^{N-2-\alpha}\left(Q_i\right)
\end{equation}
\end{thm}
Now let us show Thomae formula. As in [Na] we write the logarithmic derivative of the theta function on the divisor: $e_{\overline{\beta}}=u(\beta_iP_i)+K_{z_0}-u(\sum_{i=1}^{N}\infty_i)$
$$\frac{\partial\log\theta\left[e_{\overline{\beta}}\right]}{\partial \lambda_i}\left(0,\tau\right)={\frac{d}{dt}}\log\theta_t\left[e_{\overline{\beta}}\right](0)|_{t=0}\left(0,\tau\right)
$$
By the chain rule the last expression is:
$$
{\frac{d}{dt}}\log\theta_t\left[e_{\overline{\beta}}\right](0)|_{t=0}\left(0,\tau\right)=\sum_{1\leq k,r\leq g}\frac{\partial\log\theta\left[e_{\overline{\beta}}\right]}{\partial \tau_{kr}}\left(0\right)\frac{d\tau_{kr}}{dt}
$$
Now use the heat equation to rewrite the last expression as:
$$
\frac{1}{2}\sum_{1\leq k,r \leq g}\frac{1}{\theta\left[e_{\overline{\beta}}\right]\left(0,\tau\right)}\frac{\partial\theta\left[e_{\overline{\beta}}\right]}{\partial z_k \partial z_r}\left(0\right)\frac{d\tau_{kr}}{dt}
$$
We showed in (\textbf{~\ref{thetadervanish}})
$$\frac{\partial \theta[e]_{\overline{\beta}}}{\partial z_i}(0)=0$$ The last sum equals:
$$\frac{1}{2}\sum_{1\leq k,r\leq g}\frac{\partial\log\theta\left[e_{\overline{\beta}}\right]}{\partial z_k \partial z_r}\left(0\right)\frac{d\tau_{kr}}{dt}.$$
Use the theorem \textbf{~\ref{variational}}, and corollaries \textbf{~\ref{CorYaacov}}, \textbf{~\ref{atomic}} to conclude that:
$$\frac{\partial\log\theta\left[e_{\overline{\beta}}\right]}{\partial \lambda_i}\left(0,\tau\right)=\frac{1}{2}\frac{\partial}{\partial \lambda_i}\log \det C+\sum_{j=1,j\neq i}^{m}\frac{q(\beta_i,\beta_j)}{\lambda_i-\lambda_j}+\frac{1}{2}\sum_{j=1,j\neq i}^{j=m}
\frac{\gamma_{ij}}{\lambda_i-\lambda_j}$$
Integrate the system of first order differential equations to get the following theorem:
\begin{thm} \label{three}
Let $\beta_i$ be integer numbers and $0\leq\beta_i\leq N-1.$ such that $$\sum_{i=0}^{m}\overline{\beta_i+kR_i}=\frac{r}{2} .$$ Then there is a complex number $\alpha$ such that:
\begin{equation}
\theta\left[u(\sum_{j=1}^m \beta_iP_i)+K_{z_0}-u(\sum_{i=1}^N\infty_i)\right]=
\alpha\sqrt{\det
C}\times{\prod_{i,j=1..m,i\neq j}(\lambda_i-\lambda_j)}^{\beta_{ij}+\frac{\gamma_{ij}}{2}}
\end{equation}
where $\beta_{lj}=\sum_{k=0}^{N-1}\left(\left\{\beta_l+kR_l\right\}-\frac{N-1}{2N}\right)\left(\left\{\beta_j+kR_j\right\}-\frac{N-1}{2N}\right)$
and $\gamma_{lj}=\sum_{w=0}^{N-1}\left\{wR_l\right\}\left\{wR_j\right\}.$
\end{thm}

\end{document}